\documentclass[11pt]{amsart}
\usepackage{latexsym}
\usepackage{amsmath}
\usepackage{amssymb}

\newcommand{\eproof}{\mbox{\ }\hfill $\Box$ \par \vskip 10pt}

\newtheorem{Theorem}{Theorem}[section]
\newtheorem{lemma}[Theorem]{Lemma}

\newtheorem{corol}[Theorem]{Corollary}

\numberwithin{equation}{section}

\def\cal{\mathcal}

\begin{document}

\title[Semi-classical resolvent estimates]{Semi-classical resolvent estimates for short-range $L^\infty$ potentials. II}

\author[G. Vodev]{Georgi Vodev}

\address {Universit\'e de Nantes, Laboratoire de Math\'ematiques Jean Leray, 2 rue de la Houssini\`ere, BP 92208, 44322 Nantes Cedex 03, France}
\email{Georgi.Vodev@univ-nantes.fr}

\date{}

\begin{abstract} We prove semi-classical resolvent estimates for real-valued potentials 
$V\in L^\infty(\mathbb{R}^n)$, $n\ge 3$, of the form $V=V_L+V_S$, where $V_L$ is a long-range potential
which is $C^1$ with respect to the radial variable, while $V_S$ is a short-range potential 
satisfying $V_S(x)={\cal O}\left(\langle x\rangle^{-\delta}\right)$ with $\delta>1$. 

\quad

Key words: Schr\"odinger operator, resolvent estimates, short-range potentials.
\end{abstract} 

\maketitle

\setcounter{section}{0}
\section{Introduction and statement of results}

The goal of this paper is to extend the semi-classical resolvent estimates obtained recently in \cite{kn:KV}, \cite{kn:S2} 
and \cite{kn:V2} to a larger class of potentials. We are going to 
study the resolvent of the Schr\"odinger operator
$$P(h)=-h^2\Delta+V(x)$$
where $0<h\ll 1$ is a semi-classical parameter, $\Delta$ is the negative Laplacian in 
$\mathbb{R}^n$, $n\ge 3$, and $V\in L^\infty(\mathbb{R}^n)$ is a real-valued potential 
 of the form $V=V_L+V_S$, where $V_L\in C^1([r_0,+\infty))$ with respect to the radial variable $r=|x|$,
 $r_0>0$ being some constant, is a long-range potential, while $V_S$ is a short-range potential 
satisfying
\begin{equation}\label{eq:1.1}
|V_S(x)|\le C_1(|x|+1)^{-\delta}
\end{equation}
with some constants $C_1>0$ and $\delta>1$. We suppose that there exists a decreasing function $p(r)>0$,
$p(r)\to 0$ as $r\to\infty$, such that
\begin{equation}\label{eq:1.2}
V_L(x)\le p(|x|)\quad\mbox{for}\quad |x|\ge r_0.
\end{equation}
We also suppose that 
\begin{equation}\label{eq:1.3}
\partial_rV_L(x)\le  C_2(|x|+1)^{-\beta}\quad\mbox{for}\quad |x|\ge r_0
\end{equation}
with some constants $C_2>0$ and $\beta>1$.
As in \cite{kn:V2} we introduce the quantity
$$g_s^\pm(h,\theta):=\log\left\|(|x|+1)^{-s}(P(h)-E\pm i\theta)^{-1}(|x|+1)^{-s}
\right\|_{L^2\to L^2}$$
where $L^2:=L^2(\mathbb{R}^n)$, $0<\theta<1$, $s>1/2$ is independent of $h$ and $E>0$ is a fixed energy level independent of $h$.
Our first result is the following

\begin{Theorem} Suppose the conditions (\ref{eq:1.1}), (\ref{eq:1.2}) and (\ref{eq:1.3}) fulfilled with $\delta$ and $\beta$
satisfying the condition
\begin{equation}\label{eq:1.4}
\delta>3,\quad \beta\ge 3.
\end{equation}
Then there exist constants $C>0$ and $h_0>0$ independent of $h$ and $\theta$ but depending on $s$ and $E$ such 
that for all $0<h\le h_0$ we have the bound
\begin{equation}\label{eq:1.5}
g_s^\pm(h,\theta)\le Ch^{-4/3}\log(h^{-1}).
\end{equation}
\end{Theorem}

When $V_S\equiv 0$ and $V_L$ satisfying conditions similar to (\ref{eq:1.2}) and (\ref{eq:1.3}), it is proved in
\cite{kn:D} when $n\ge 3$ and in \cite{kn:S1} when $n=2$ that
\begin{equation}\label{eq:1.6}
g_s^\pm(h,\theta)\le Ch^{-1}
\end{equation}
with some constant $C>0$ independent of $h$ and $\theta$. Previously, the bound (\ref{eq:1.6}) 
was proved for smooth potentials
in \cite{kn:B2} and an analog of (\ref{eq:1.6}) for H\"older potentials was proved in \cite{kn:V1}. 
A high-frequency analog of (\ref{eq:1.6}) on Riemannian manifolds was also proved in 
\cite{kn:B1} and \cite{kn:CV}.  When $V_L\equiv 0$ and $V_S$ satisfying the condition (\ref{eq:1.1}) with $\delta>3$, 
the bound (\ref{eq:1.5}) has been
recently proved in \cite{kn:V2}. Previously, (\ref{eq:1.5}) was proved in \cite{kn:KV} and \cite{kn:S2} 
for real-valued compactly supported $L^\infty$ potentials.
When $n=1$ it was shown in \cite{kn:DZ} that we have the better bound (\ref{eq:1.6}) instead of (\ref{eq:1.5}). 
The method we use to prove Theorem 1.1 also allows us to get resolvent bounds when the condition (\ref{eq:1.4}) is not satisfied,
which however are much weaker than the bound (\ref{eq:1.5}). More precisely, we have the following

\begin{Theorem} Suppose the conditions (\ref{eq:1.1}), (\ref{eq:1.2}) and (\ref{eq:1.3}) fulfilled 
with $\delta$ and $\beta$
satisfying either the condition
\begin{equation}\label{eq:1.7}
1<\delta\le 3,\quad \beta>1,
\end{equation}
or the condition
\begin{equation}\label{eq:1.8}
\delta>3,\quad 1<\beta<3.
\end{equation}
Then,  
there exist constants $C>0$ and $h_0>0$ independent of $h$ and $\theta$ but depending on $s$ and $E$ such 
that for all $0<h\le h_0$ we have the bounds 
\begin{equation}\label{eq:1.9}
g_s^\pm(h,\theta)\le 
\left\{
\begin{array}{ll}
 Ch^{-\frac{2}{3}-m_1}& \mbox{if (\ref{eq:1.7}) holds},\\
 Ch^{-\frac{4}{3}-\frac{1}{2}(3-\beta)m_2}& \mbox{if (\ref{eq:1.8}) holds},\\
\end{array}
\right.
\end{equation}
 where
 $$m_1=\max\left\{\frac{7}{3(\delta-1)},\frac{4}{3(\beta-1)}\right\}\ge \frac{7}{6}$$
 and
 $$m_2=\max\left\{\frac{1}{\delta-\beta},\frac{4}{3(\beta-1)}\right\}> \frac{2}{3}.$$
\end{Theorem}

Clearly, this theorem implies the following

\begin{corol} Suppose that $V_L\equiv 0$ and let $V=V_S$ satisfy 
the condition (\ref{eq:1.1}) with $1<\delta\le 3$. 
Then,  
there exist constants $C>0$ and $h_0>0$ independent of $h$ and $\theta$ but depending on $s$ and $E$ such 
that for all $0<h\le h_0$ we have the bound 
\begin{equation}\label{eq:1.10}
g_s^\pm(h,\theta)\le 
 Ch^{-\frac{2\delta+5}{3(\delta-1)}}.
\end{equation}
\end{corol}

To prove the above theorems we follow the same strategy as in \cite{kn:V2} which in turn is inspired by the paper \cite{kn:S2}. 
It consists of using Carleman estimates with phase and weight functions, denoted by $\varphi$ and $\mu$ below, depending only on
the radial variable $r$ and the parameter $h$, which have very weak regularity. It turns out that it suffices to choose $\varphi$ belonging only to 
$C^1$ and $\mu$ only continuous. Thus we get derivatives $\varphi''$ and $\mu'$ belonging to $L^\infty$, which proves
sufficient for the Carleman estimates to hold. Note that higher derivatives of $\varphi$ and $\mu$ are not involved in the proof of 
the Carleman estimates (see the proof of Theorem 3.1 below). In order to be able to prove the Carleman estimates the functions 
$\varphi$ and $\mu$ must satisfy some conditions (see the inequalities (\ref{eq:2.2}) and (\ref{eq:2.8}) below). On the other hand, 
to get as good resolvent bounds as possible we are looking for a phase function $\varphi$ such that $\max\varphi$ is as small as possible. 
The construction of such phase and weight functions is carried out in Section 2 following that one in \cite{kn:V2}. However, here 
the construction is more complicated due to the more general class the potential belongs to. 

It is not clear if the bounds (\ref{eq:1.5}) and (\ref{eq:1.9}) are optimal for $L^\infty$ potentials. In any case, they seem hard to improve
unless one menages to construct a better phase function. By contrast, the optimality of the bound (\ref{eq:1.6}) for smooth potentials 
is well known (e.g. see \cite{kn:DDZ}).

\section{The construction of the phase and weight functions revisited} 

We will follow closely the construction in Section 2 of \cite{kn:V2} making some suitable modifications in order to adapte it
to the more general class of potentials we consider in the present paper. 
We will first construct the weight function $\mu$ as follows:
$$\mu(r)=
\left\{
\begin{array}{lll}
 (r+1)^{2k}-1&\mbox{for}& 0\le r\le a,\\
 (a+1)^{2k}-1+(a+1)^{-2s+1}-(r+1)^{-2s+1}&\mbox{for}& r\ge a,
\end{array}
\right.
$$
where $a=h^{-m}$ with
$$m=\left\{
\begin{array}{lll}
 m_0+\epsilon T_0& \mbox{if (\ref{eq:1.4}) holds},\\
 m_1+\epsilon T_1& \mbox{if (\ref{eq:1.7}) holds},\\
 m_2+\epsilon T_2& \mbox{if (\ref{eq:1.8}) holds},\\
\end{array}
\right.
$$
where $\epsilon=\left(\log\frac{1}{h}\right)^{-1}$, $m_0=\max\left\{\frac{2}{3},\frac{1}{\delta-3}\right\}$, $m_1$ and $m_2$
are as in Theorem 1.2 and $T_0, T_1, T_2>0$ are parameters independent of $h$ to be fixed in the proof of Lemma 2.3 below. Furthermore,  
$$k=\left\{
\begin{array}{lll}
 1-\epsilon t_0& \mbox{if (\ref{eq:1.4}) holds},\\
 \frac{2}{3m_1}-\epsilon t_1& \mbox{if (\ref{eq:1.7}) holds},\\
  \frac{1}{2}(\beta-1)-\epsilon t_2& \mbox{if (\ref{eq:1.8}) holds},\\
\end{array}
\right.
$$
and
\begin{equation}\label{eq:2.1}
s=\frac{1+\epsilon}{2}
\end{equation}
where $t_0, t_1, t_2>1$ are parameters independent of $h$ to be fixed in the proof of Lemma 2.3 below. 

Clearly, the first derivative (in sense of distributions) of $\mu$ satisfies
$$\mu'(r)=
\left\{
\begin{array}{lll}
 2k(r+1)^{2k-1}&\mbox{for}& 0\le r<a,\\
 (2s-1)(r+1)^{-2s}&\mbox{for}& r>a.
\end{array}
\right.
$$
The following properties of the functions $\mu$ and $\mu'$ are essential to prove the Carleman estimates in the next section.

\begin{lemma} For all $r>0$, $r\neq a$, we have the inequalities
\begin{equation}\label{eq:2.2}
2r^{-1}\mu(r)-\mu'(r)\ge 0,
\end{equation}
\begin{equation}\label{eq:2.3}
\mu'(r)\ge \epsilon(r+1)^{-2s},
\end{equation}
\begin{equation}\label{eq:2.4}
\frac{\mu(r)^2}{\mu'(r)}\le 2\epsilon^{-1}a^{4k}(r+1)^{2s}.
\end{equation}
\end{lemma}

{\it Proof.} It is easy to see that for $r<a$ (\ref{eq:2.2}) follows from the inequality
$$f(r):=1+(1-k)r-(r+1)^{1-2k}\ge 0$$
for all $r\ge 0$ and $0\le k\le 1$. It is obvious for $1/2\le k\le 1$, while for $0\le k<1/2$ we have
$$f'(r)=1-k-(1-2k)(r+1)^{-2k}\ge k\ge 0.$$
Hence in this case the function $f$ is increasing, which implies $f(r)\ge f(0)=0$ as desired.

For $r>a$ the left-hand side of (\ref{eq:2.2}) is bounded
from below by
$$ 2r^{-1}((a+1)^{2k}-1-s)>0$$
provided $a$ is taken large enough. Furthermore, we clearly have $\mu'(r)\ge 2k(r+1)^{-1}$ for $r<a$, and hence (\ref{eq:2.3}) holds in this case,
provided $\epsilon$ is taken small enough.  
For $r>a$ the bound (\ref{eq:2.3}) is trivial.   
The bound (\ref{eq:2.4}) follows from (\ref{eq:2.3}) and the observation that $\mu(r)^2\le (a+1)^{4k}\le 2a^{4k}$ for all $r$. 
\eproof

We now turn to the 
construction of the phase function
$\varphi\in C^1([0,+\infty))$ such that $\varphi(0)=0$ and $\varphi(r)>0$ for $r>0$. 
We define the first derivative of $\varphi$ by
$$\varphi'(r)=
\left\{
\begin{array}{lll}
 \tau(r+1)^{-k}- \tau(a+1)^{-k}&\mbox{for}& 0\le r\le a,\\
 0&\mbox{for}& r\ge a,
\end{array}
\right.
$$
where 
\begin{equation}\label{eq:2.5}
\tau=\tau_0h^{-1/3}
\end{equation}
with some parameter $\tau_0\gg 1$ independent of $h$ to be fixed in Lemma 2.3 below.
Clearly, the first derivative of $\varphi'$ satisfies
$$\varphi''(r)=
\left\{
\begin{array}{lll}
 -k\tau(r+1)^{-k-1}&\mbox{for}& 0\le r<a,\\
 0&\mbox{for}& r>a.
\end{array}
\right.
$$

\begin{lemma} For all $r\ge 0$ we have the bounds
\begin{equation}\label{eq:2.6}
h^{-1}\varphi(r)\lesssim 
\left\{
\begin{array}{lll}
 h^{-4/3}\log\frac{1}{h}& \mbox{if (\ref{eq:1.4}) holds},\\
 h^{-\frac{2}{3}-m_1}& \mbox{if (\ref{eq:1.7}) holds},\\
 h^{-\frac{4}{3}-\frac{1}{2}(3-\beta)m_2}& \mbox{if (\ref{eq:1.8}) holds}.\\
\end{array}
\right.
\end{equation}
\end{lemma}

{\it Proof.} Since $k<1$ we have
$$\max\varphi=\int_0^a\varphi'(r)dr\le \tau\int_0^a (r+1)^{-k}dr\le
 \frac{\tau}{1-k}(a+1)^{1-k}$$
 $$\lesssim 
\left\{
\begin{array}{ll}
 \tau\epsilon^{-1}& \mbox{if (\ref{eq:1.4}) holds},\\
 \tau a^{1-k}& \mbox{otherwise},\\
\end{array}
\right.
$$
where we have used that $a^\epsilon={\cal O}(1)$ and $(1-k)^{-1}={\cal O}(\epsilon^{-1})$ if (\ref{eq:1.4}) holds, 
$(1-k)^{-1}={\cal O}(1)$ in the other two cases. 
In view of the choice of $\epsilon$, $\tau$ and $a$, we get the bounds
\begin{equation}\label{eq:2.7}
h^{-1}\varphi(r)\lesssim 
\left\{
\begin{array}{ll}
 h^{-4/3}\log\frac{1}{h}& \mbox{if (\ref{eq:1.4}) holds},\\
  h^{-4/3-m(1-k)}& \mbox{otherwise}.\\
\end{array}
\right.
\end{equation}
Since 
$$m(1-k)=\left\{
\begin{array}{ll}
 m_1-\frac{2}{3}+{\cal O}(\epsilon)& \mbox{if (\ref{eq:1.7}) holds},\\
 \frac{1}{2}(3-\beta)m_2+{\cal O}(\epsilon)& \mbox{if (\ref{eq:1.8}) holds},\\
\end{array}
\right.
$$
(\ref{eq:2.7}) clearly implies (\ref{eq:2.6}).
\eproof

Let $\phi\in C_0^\infty([1,2])$, $\phi\ge 0$, be a real-valued function independent of $h$ such that $\int_{-\infty}^\infty\phi(\sigma)d\sigma=1$.
Given a parameter $b\gg r_0$ to be fixed in the proof of Theorem 3.1 below, independent of $h$, set
$$\psi_b(r)=b^{-1}\int_{r}^\infty\phi(\sigma/b)d\sigma.$$
Clearly, we have $0\le\psi_b\le 1$ and $\psi_b(r)=1$ for $r\le b$, $\psi_b(r)=0$ for $r\ge 2b$.
For $r>0$, $r\neq a$, set
$$A(r)=\left(\mu\varphi'^2\right)'(r)$$
and
$$B(r)=\frac{3\left(\mu(r)\left(h^{-1}C_1(r+1)^{-\delta}+h^{-1}Q_b\psi_b(r)+|\varphi''(r)|\right)\right)^2}{h^{-1}\varphi'(r)\mu(r)+\mu'(r)}$$
$$+\mu(r)(1-\psi_b(r))C_2(r+1)^{-\beta}$$
where $Q_b\ge 0$ is some constant depending only on $b$. 
The following lemma will play a crucial role in the proof of the Carleman estimates in the next section.

\begin{lemma} There exist constants $b_0=b_0(E)>0$, $\tau_0=\tau_0(b,E)>0$
 and 
$h_0=h_0(b,E)>0$ so that for $\tau$ satisfying (\ref{eq:2.5}) and for all $b\ge b_0$, $0<h\le h_0$ we have the inequality
\begin{equation}\label{eq:2.8}
A(r)-B(r)\ge -\frac{E}{2}\mu'(r)
\end{equation}
for all $r>0$, $r\neq a$.
\end{lemma}

{\it Proof.} For $r<a$ we have
$$A(r)=-\left(\varphi'^2\right)'(r)+\tau^2\partial_r\left(1-(r+1)^k(a+1)^{-k}\right)^2$$ 
$$=-2\varphi'(r)\varphi''(r)-2k\tau^2(r+1)^{k-1}(a+1)^{-k}\left(1-(r+1)^k(a+1)^{-k}\right)$$
$$\ge 2k\tau(r+1)^{-k-1}\varphi'(r)-2k\tau^2(r+1)^{k-1}(a+1)^{-k}$$
$$\ge 2k\tau(r+1)^{-k-1}\varphi'(r)-\tau^2 a^{-k}\mu'(r).$$
Taking into account the definition of the parameters $a$ and $\tau$ we conclude
\begin{equation}\label{eq:2.9}
A(r)\ge 2k\tau(r+1)^{-k-1}\varphi'(r)-{\cal O}(h^{km-2/3})\mu'(r)
\end{equation}
for all $r<a$. Observe now that if (\ref{eq:1.4}) holds, we have
$$km-2/3=m_0-2/3+\epsilon(T_0-m_0t_0)-{\cal O}(\epsilon^2)\ge \epsilon m_0t_0$$
provided we take $T_0=3m_0t_0$ and $\epsilon$ small enough. 
If (\ref{eq:1.7}) holds, we have
$$km-2/3=\frac{2\epsilon T_1}{3m_1}-\epsilon m_1t_1-{\cal O}(\epsilon^2).$$
We take now $T_1=6m_1^2t_1$. Then
$$km-2/3=3\epsilon m_1t_1-{\cal O}(\epsilon^2)\ge \epsilon m_1t_1$$
provided $\epsilon$ is taken small enough. 
On the other hand, if (\ref{eq:1.8}) holds, we have
$$km-2/3=\frac{(\beta-1)}{2}m_2-\frac{2}{3}+\frac{(\beta-1)}{2}\epsilon T_2-\epsilon m_2t_2
-{\cal O}(\epsilon^2)$$
$$\ge 2\epsilon m_2t_2-{\cal O}(\epsilon^2)\ge \epsilon m_2t_2$$
provided we take
$$T_2=\frac{6m_2t_2}{\beta-1}$$
and $\epsilon$ small enough. Using that $h^{\epsilon t}=e^{-t}$ we conclude 
\begin{equation}\label{eq:2.10}
h^{km-2/3}\le
\left\{
\begin{array}{lll}
e^{-t_0m_0} & \mbox{if (\ref{eq:1.4}) holds} ,\\
e^{-t_1m_1}& \mbox{if (\ref{eq:1.7}) holds},\\
e^{-t_2m_2}& \mbox{if (\ref{eq:1.8}) holds}.\\
\end{array}
\right.
\end{equation}
Taking $t_0$, $t_1$ and $t_2$ large enough, independent of $h$, we obtain from (\ref{eq:2.9}) and (\ref{eq:2.10}) that in all cases
we have the estimate
\begin{equation}\label{eq:2.11}
A(r)\ge 2k\tau(r+1)^{-k-1}\varphi'(r)-\frac{E}{4}\mu'(r)
\end{equation}
for all $r<a$.
We will now bound the function $B$ from above. Note that taking $h$ small enough we can arrange that $2b<a/2$. 
Let first $0<r\le \frac{a}{2}$. Since in this case we have
$$\varphi'(r)\ge \widetilde C\tau(r+1)^{-k}$$
with some constant $\widetilde C>0$, we obtain
$$B(r)\lesssim \frac{\mu(r)\left(h^{-2}\widetilde Q_b(r+1)^{-2\delta}+\varphi''(r)^2\right)}{h^{-1}\varphi'(r)}+\mu(r)(1-\psi_b(r))(r+1)^{-\beta}$$ 
 $$\lesssim \widetilde Q_b(\tau h)^{-1}\frac{\mu(r)(r+1)^{1+k-2\delta}}{\varphi'(r)^2}\tau(r+1)^{-k-1}\varphi'(r)
 + h\frac{\mu(r)\varphi''(r)^2}{\mu'(r)\varphi'(r)}\mu'(r)$$ 
 $$+(1-\psi_b(r))(r+1)^{2k-\beta}$$ 
$$\lesssim \widetilde Q_b\tau^{-3}h^{-1}(r+1)^{1+5k-2\delta}\tau(r+1)^{-k-1}\varphi'(r)+\tau h\mu'(r)$$
$$+(1-\psi_b(r))(r+1)^{1-\beta}\mu'(r)$$ 
$$\lesssim \widetilde Q_b\tau_0^{-3}\tau(r+1)^{-k-1}\varphi'(r)+ (\tau_0h^{2/3}+b^{-\beta+1})\mu'(r)$$
where $\widetilde Q_b>0$ is some constant depending only on $b$ and we have used that $k<(2\delta-1)/5$ in all three cases. 
Taking $h$ small enough, depending on $\tau_0$, and $b$ big enough, independent of $h$ and $\tau_0$, we get the bound
\begin{equation}\label{eq:2.12}
B(r)\le C\widetilde Q_b\tau_0^{-3}\tau(r+1)^{-k-1}\varphi'(r)+\frac{E}{4}\mu'(r)
\end{equation}
with some constant $C>0$. 
In this case we get (\ref{eq:2.8}) from (\ref{eq:2.11}) and (\ref{eq:2.12}) by taking $\tau_0$ big enough
depending on $b$ and $C$ but independent of $h$.

Let now $\frac{a}{2}<r<a$. Then we have the bound
$$B(r)\lesssim \left(\frac{\mu(r)}{\mu'(r)}\right)^2\left(h^{-1}(r+1)^{-\delta}+|\varphi''(r)|\right)^2\mu'(r)+(r+1)^{-\beta+1}\mu'(r)$$
$$\lesssim \left(h^{-2}(r+1)^{2-2\delta}+\tau^2(r+1)^{-2k}\right)\mu'(r)+a^{-\beta+1}\mu'(r)$$
$$\lesssim \left(h^{-2}a^{2-2\delta}+\tau^2a^{-2k}\right)\mu'(r)+a^{-\beta+1}\mu'(r)$$
$$\lesssim \left(h^{2m(\delta-1)-2}+h^{2km-2/3}+h^{m(\beta-1)}\right)\mu'(r)\le \frac{E}{4}\mu'(r)$$
provided $h$ is taken small enough. Again, this bound together with (\ref{eq:2.11}) imply (\ref{eq:2.8}).

It remains to consider the case $r>a$. Using that $\mu={\cal O}(a^{2k})$ together with (\ref{eq:2.3}) we get
$$B(r)\lesssim \frac{\left(\mu(r)\left(h^{-1}(r+1)^{-\delta}\right)\right)^2}{\mu'(r)}+(r+1)^{-\beta}\mu(r)$$
$$\lesssim h^{-2}a^{4k}(r+1)^{4s-2\delta}\mu'(r)+a^{2k}(r+1)^{2s-\beta}\mu'(r)$$
$$\lesssim \left(h^{-2}a^{4k+4s-2\delta}+a^{2k+2s-\beta}\right)\mu'(r)$$
$$\lesssim \left(h^{2m(\delta-2k-2s)-2}+h^{m(\beta-2k-2s)}\right)\mu'(r).$$
When (\ref{eq:1.4}) holds we have 
$$2k+2s=3-(2t_0-1)\epsilon<3-t_0\epsilon$$
 and hence
$$m(\delta-2k-2s)-1\ge m_0(\delta-3)-1+\epsilon m_0t_0\ge \epsilon m_0t_0$$
and
$$m(\beta-2k-2s)\ge m_0(\beta-3+\epsilon t_0)\ge \epsilon m_0t_0.$$
When (\ref{eq:1.7}) holds we have 
$$2k+2s=\frac{4}{3m_1}+1-(2t_1-1)\epsilon<\frac{4}{3m_1}+1-t_1\epsilon$$
 and hence
$$m(\delta-2k-2s)-1\ge m_1(\delta-1-\frac{4}{3m_1}+\epsilon t_1)-1= m_1(\delta-1)-\frac{7}{3}+\epsilon m_1t_1\ge \epsilon m_1t_1.$$
 In this case we also have
$$m(\beta-2k-2s)\ge m_1(\beta-1-\frac{4}{3m_1}+\epsilon t_1)
= m_1(\beta-1)-\frac{4}{3}+\epsilon m_1t_1\ge \epsilon m_1t_1.$$
When (\ref{eq:1.8}) holds we have
$$2k+2s=\beta-(2t_2-1)\epsilon<\beta-t_2\epsilon$$
 and hence
$$m(\delta-2k-2s)-1\ge m_2(\delta-\beta+\epsilon t_2)-1\ge \epsilon m_2t_2$$
and
$$m(\beta-2k-2s)\ge \epsilon m_2t_2.$$
We conclude from the above inequalities that
\begin{equation}\label{eq:2.13}
h^{2m(\delta-2k-2s)-2}+h^{m(\beta-2k-2s)}\le
\left\{
\begin{array}{lll}
2e^{-t_0m_0}& \mbox{if (\ref{eq:1.4}) holds} ,\\
2e^{-t_1m_1}& \mbox{if (\ref{eq:1.7}) holds},\\
2e^{-t_2m_2}& \mbox{if (\ref{eq:1.8}) holds}.\\
\end{array}
\right.
\end{equation}
It follows from (\ref{eq:2.13}) that taking $t_0$, $t_1$ and $t_2$ large enough, independent of $h$, we can arrange the bound
\begin{equation}\label{eq:2.14}
B(r)\le \frac{E}{2}\mu'(r).
\end{equation}
 Since in this case $A(r)=0$, the bound (\ref{eq:2.14}) clearly implies (\ref{eq:2.8}).
\eproof

\section{Carleman estimates} 

In this section we will prove the following

\begin{Theorem} Suppose (\ref{eq:1.1}), (\ref{eq:1.2}) and (\ref{eq:1.3}) fulfilled and let $s$ satisfy (\ref{eq:2.1}). Then, for all functions
$f\in H^2(\mathbb{R}^n)$ such that $(|x|+1)^{s}(P(h)-E\pm i\theta)f\in L^2$ and for all
$0<h\le h_0$, $0<\theta\le \epsilon ha^{-2k}$, we have the estimate 
 \begin{equation}\label{eq:3.1}
\|(|x|+1)^{-s}e^{\varphi/h}f\|_{L^2}\le Ca^{2k}(\epsilon h)^{-1}\|(|x|+1)^{s}e^{\varphi/h}(P(h)-E\pm i\theta)f\|_{L^2}$$ $$
+Ca^k\tau\left(\frac{\theta}{\epsilon h}\right)^{1/2}\|e^{\varphi/h}f\|_{L^2}
\end{equation}
with a constant $C>0$ independent of $h$, $\theta$ and $f$.
\end{Theorem}

{\it Proof.} We will adapt the proof of Theorem 3.1 of \cite{kn:V2} to this more general case. 
We pass to the polar coordinates $(r,w)\in\mathbb{R}^+\times\mathbb{S}^{n-1}$, $r=|x|$, $w=x/|x|$, and
recall that $L^2(\mathbb{R}^n)=L^2(\mathbb{R}^+\times\mathbb{S}^{n-1}, r^{n-1}drdw)$. In what follows we denote by $\|\cdot\|$ and $\langle\cdot,\cdot\rangle$
the norm and the scalar product in $L^2(\mathbb{S}^{n-1})$. We will make use of the identity
\begin{equation}\label{eq:3.2}
 r^{(n-1)/2}\Delta  r^{-(n-1)/2}=\partial_r^2+\frac{\widetilde\Delta_w}{r^2}
\end{equation}
where $\widetilde\Delta_w=\Delta_w-\frac{1}{4}(n-1)(n-3)$ and $\Delta_w$ denotes the negative Laplace-Beltrami operator
on $\mathbb{S}^{n-1}$. Set $u=r^{(n-1)/2}e^{\varphi/h}f$ and
$${\cal P}^\pm(h)=r^{(n-1)/2}(P(h)-E\pm i\theta) r^{-(n-1)/2},$$
$${\cal P}^\pm_\varphi(h)=e^{\varphi/h}{\cal P}^\pm(h)e^{-\varphi/h}.$$
Using (\ref{eq:3.2}) we can write the operator ${\cal P}^\pm(h)$ in the coordinates $(r,w)$ as follows
$${\cal P}^\pm(h)={\cal D}_r^2+\frac{\Lambda_w}{r^2}-E\pm i\theta +V$$
where we have put ${\cal D}_r=-ih\partial_r$ and $\Lambda_w=-h^2\widetilde\Delta_w$. Since the function $\varphi$
depends only on the variable $r$, this implies
$${\cal P}^\pm_\varphi(h)={\cal D}_r^2+\frac{\Lambda_w}{r^2}-E\pm i\theta -\varphi'^2+h\varphi''+
2i\varphi'{\cal D}_r+V.$$
We now write $V=\widetilde V_S+\widetilde V_L$ with 
$$\widetilde V_S(x)=V_S(x)+\psi_b(|x|)V_L(x)$$
 and $$\widetilde V_L(x)=(1-\psi_b(|x|))V_L(x).$$ 
For $r>0$, $r\neq a$, introduce the function
$$F(r)=-\langle (r^{-2}\Lambda_w-E-\varphi'(r)^2+\widetilde V_L(r,\cdot))u(r,\cdot),u(r,\cdot)\rangle+\|{\cal D}_ru(r,\cdot)\|^2$$
where $\widetilde V_L(r,w):=\widetilde V_L(rw)$. 
It is easy to check that its first derivative is given by
$$F'(r)=\frac{2}{r}\langle r^{-2}\Lambda_wu(r,\cdot),u(r,\cdot)\rangle
+((\varphi')^2-\widetilde V_L)'\|u(r,\cdot)\|^2$$
$$-2h^{-1}{\rm Im}\,\langle {\cal P}^\pm_\varphi(h)u(r,\cdot),{\cal D}_ru(r,\cdot)\rangle$$
$$\pm 2\theta h^{-1}{\rm Re}\,\langle u(r,\cdot),{\cal D}_ru(r,\cdot)\rangle+4h^{-1}\varphi'\|{\cal D}_ru(r,\cdot)\|^2$$ 
$$+2h^{-1}{\rm Im}\,\langle (\widetilde V_S+h\varphi'')u(r,\cdot),{\cal D}_ru(r,\cdot)\rangle.$$
Thus, if $\mu$ is the function defined in the previous section, we obtain the identity
$$\mu'F+\mu F'=
(2r^{-1}\mu-\mu')\langle r^{-2}\Lambda_wu(r,\cdot),u(r,\cdot)\rangle
+(E\mu'+(\mu(\varphi')^2-\mu\widetilde V_L)')\|u(r,\cdot)\|^2$$
$$-2h^{-1}\mu{\rm Im}\,\langle {\cal P}^\pm_\varphi(h)u(r,\cdot),{\cal D}_ru(r,\cdot)\rangle$$
$$\pm 2\theta h^{-1}\mu{\rm Re}\,\langle u(r,\cdot),{\cal D}_ru(r,\cdot)\rangle+(\mu'+4h^{-1}\varphi'\mu)\|{\cal D}_ru(r,\cdot)\|^2$$ 
$$+2h^{-1}\mu{\rm Im}\,\langle (\widetilde V_S+h\varphi'')u(r,\cdot),{\cal D}_ru(r,\cdot)\rangle.$$
Using that $\Lambda_w\ge 0$ together with (\ref{eq:2.2}) we get the inequality
$$\mu'F+\mu F'\ge (E\mu'+(\mu(\varphi')^2-\mu\widetilde V_L)')\|u(r,\cdot)\|^2+(\mu'+4h^{-1}\varphi'\mu)\|{\cal D}_ru(r,\cdot)\|^2$$
$$-\frac{3h^{-2}\mu^2}{\mu'}\|{\cal P}^\pm_\varphi(h)u(r,\cdot)\|^2-\frac{\mu'}{3}\|{\cal D}_ru(r,\cdot)\|^2$$
$$-\theta h^{-1}\mu\left(\|u(r,\cdot)\|^2+\|{\cal D}_ru(r,\cdot)\|^2\right)$$
$$-3h^{-2}\mu^2(\mu'+4h^{-1}\varphi'\mu)^{-1}\|(\widetilde V_S+h\varphi'')u(r,\cdot)\|^2
-\frac{1}{3}(\mu'+4h^{-1}\varphi'\mu)\|{\cal D}_ru(r,\cdot)\|^2.$$
In view of the assumptions (\ref{eq:1.2}) and (\ref{eq:1.3}) we have
$$(\mu\widetilde V_L)'=\mu'\widetilde V_L+\mu\widetilde V'_L=\mu'(1-\psi_b)V_L-\mu\psi'_bV_L+\mu(1-\psi_b)V'_L$$
$$\le \mu'(1-\psi_b)p(r)+\mu b^{-1}\phi(r/b)p(r)+\mu(1-\psi_b)C_2(r+1)^{-\beta}$$
$$\le \mu'(1-\psi_b)p(b)+{\cal O}(r) b^{-1}\phi(r/b)p(b)\mu'+\mu(1-\psi_b)C_2(r+1)^{-\beta}$$
$$\le {\cal O}(1)p(b)\mu'+\mu(1-\psi_b)C_2(r+1)^{-\beta}\le \frac{E}{3}\mu'+\mu(1-\psi_b)C_2(r+1)^{-\beta}$$
provided $b$ is taken large enough. Observe also that the assumption (\ref{eq:1.1}) yields
$$|\widetilde V_S|\le |V_S|+\psi_b|V_L|\le C_1(r+1)^{-\delta}+Q_b\psi_b$$
where $Q_b=\sup_{|x|\le 2b}|V_L(x)|$. Combining the above inequalities we get
 $$\mu'F+\mu F'\ge\left(\frac{2E}{3}\mu'+(\mu(\varphi')^2)'\right)\|u(r,\cdot)\|^2$$
 $$-\left(3\mu^2(\mu'+h^{-1}\varphi'\mu)^{-1}(h^{-1}C_1(r+1)^{-\delta}+h^{-1}Q_b\psi_b
 +|\varphi''|)^2+\mu(1-\psi_b)C_2(r+1)^{-\beta}\right)\|u(r,\cdot)\|^2$$
$$-\frac{3h^{-2}\mu^2}{\mu'}\|{\cal P}^\pm_\varphi(h)u(r,\cdot)\|^2
-\theta h^{-1}\mu\left(\|u(r,\cdot)\|^2+\|{\cal D}_ru(r,\cdot)\|^2\right)$$
$$=\left(\frac{2E}{3}\mu'+A(r)-B(r)\right)\|u(r,\cdot)\|^2$$
$$-\frac{3h^{-2}\mu^2}{\mu'}\|{\cal P}^\pm_\varphi(h)u(r,\cdot)\|^2
-\theta h^{-1}\mu\left(\|u(r,\cdot)\|^2+\|{\cal D}_ru(r,\cdot)\|^2\right).$$
 Now we use Lemma 2.3 to conclude that
$$\mu'F+\mu F'\ge \frac{E}{6}\mu'\|u(r,\cdot)\|^2-\frac{3h^{-2}\mu^2}{\mu'}\|{\cal P}^\pm_\varphi(h)u(r,\cdot)\|^2$$ $$
-\theta h^{-1}\mu\left(\|u(r,\cdot)\|^2+\|{\cal D}_ru(r,\cdot)\|^2\right).$$
We integrate this inequality with respect to $r$ and use that, since $\mu(0)=0$, we have
$$\int_0^\infty(\mu'F+\mu F')dr=0.$$
Thus we obtain the estimate
\begin{equation}\label{eq:3.3}
\frac{E}{6}\int_0^\infty\mu'\|u(r,\cdot)\|^2dr\le 3h^{-2}\int_0^\infty\frac{\mu^2}{\mu'}
\|{\cal P}^\pm_\varphi(h)u(r,\cdot)\|^2dr$$ $$
+\theta h^{-1}\int_0^\infty\mu\left(\|u(r,\cdot)\|^2+\|{\cal D}_ru(r,\cdot)\|^2\right)dr.
\end{equation}
Using that $\mu={\cal O}(a^{2k})$ together with (\ref{eq:2.3}) and (\ref{eq:2.4}) we get from (\ref{eq:3.3})
\begin{equation}\label{eq:3.4}
\int_0^\infty(r+1)^{-2s}\|u(r,\cdot)\|^2dr\le Ca^{4k}(\epsilon h)^{-2}\int_0^\infty(r+1)^{2s}\|{\cal P}^\pm_\varphi(h)u(r,\cdot)\|^2dr$$ $$
+C\theta (\epsilon h)^{-1}a^{2k}\int_0^\infty\left(\|u(r,\cdot)\|^2+\|{\cal D}_ru(r,\cdot)\|^2\right)dr
\end{equation}
with some constant $C>0$ independent of $h$ and $\theta$. On the other hand, we have the identity
$${\rm Re}\,\int_0^\infty\langle 2i\varphi'{\cal D}_ru(r,\cdot),u(r,\cdot)\rangle dr=\int_0^\infty h\varphi''\|u(r,\cdot)\|^2dr$$
and hence
$${\rm Re}\,\int_0^\infty\langle {\cal P}^\pm_\varphi(h)u(r,\cdot),u(r,\cdot)\rangle dr
=\int_0^\infty\|{\cal D}_ru(r,\cdot)\|^2dr
+\int_0^\infty \langle r^{-2}\Lambda_wu(r,\cdot),u(r,\cdot)\rangle dr$$
$$-\int_0^\infty(E+\varphi'^2)\|u(r,\cdot)\|^2dr
+\int_0^\infty\langle Vu(r,\cdot),u(r,\cdot)\rangle dr.$$
This implies
\begin{equation}\label{eq:3.5}
\int_0^\infty\|{\cal D}_ru(r,\cdot)\|^2dr\le {\cal O}(\tau^2)\int_0^\infty\|u(r,\cdot)\|^2dr$$
$$+\gamma\int_0^\infty(r+1)^{-2s}\|u(r,\cdot)\|^2dr
+\gamma^{-1}\int_0^\infty(r+1)^{2s}\|{\cal P}^\pm_\varphi(h)u(r,\cdot)\|^2dr
\end{equation}
for every $\gamma>0$. We take now $\gamma$ small enough, independent of $h$, and recall that $\theta (\epsilon h)^{-1}a^{2k}\le 1$. Thus, combining
the estimates (\ref{eq:3.4}) and (\ref{eq:3.5}), we get
\begin{equation}\label{eq:3.6}
\int_0^\infty(r+1)^{-2s}\|u(r,\cdot)\|^2dr\le Ca^{4k}(\epsilon h)^{-2}\int_0^\infty(r+1)^{2s}\|{\cal P}^\pm_\varphi(h)u(r,\cdot)\|^2dr$$ $$
+C\theta (\epsilon h)^{-1}a^{2k}\tau^2\int_0^\infty\|u(r,\cdot)\|^2dr
\end{equation}
with a new constant $C>0$ independent of $h$ and $\theta$. Clearly, the estimate 
(\ref{eq:3.6}) implies (\ref{eq:3.1}).
\eproof

\section{Resolvent estimates}

The bounds (\ref{eq:1.5}) and (\ref{eq:1.9}) can be derived from Theorem 3.1 in the same way as in Section 4 of \cite{kn:V2}. 
Here we will sketch the proof for the sake of completeness. Observe that it follows from the estimate 
(\ref{eq:3.1}) and Lemma 2.2 that for $0<h\ll 1$, $0<\theta\le \epsilon ha^{-2k}$ and $s$ satisfying (\ref{eq:2.1}) we have the estimate 
\begin{equation}\label{eq:4.1}
\|(|x|+1)^{-s}f\|_{L^2}\le M\|(|x|+1)^{s}(P(h)-E\pm i\theta)f\|_{L^2}
+M\theta^{1/2}\|f\|_{L^2}
\end{equation}
where
$$M=\left\{
\begin{array}{lll}
 \exp\left(Ch^{-4/3}\log\frac{1}{h}\right)& \mbox{if (\ref{eq:1.4}) holds},\\
 \exp\left(Ch^{-\frac{2}{3}-m_1}\right)& \mbox{if (\ref{eq:1.7}) holds},\\
 \exp\left(Ch^{-\frac{4}{3}-\frac{1}{2}(3-\beta)m_2}\right)& \mbox{if (\ref{eq:1.8}) holds},\\
\end{array}
\right.$$
with a constant $C>0$ independent of $h$ and $\theta$. On the other hand, since the operator $P(h)$ is symmetric, we have
\begin{equation}\label{eq:4.2}
\theta\|f\|^2_{L^2}=\pm{\rm Im}\,\langle (P(h)-E\pm i\theta)f,f\rangle_{L^2}$$
$$\le (2M)^{-2}\|(|x|+1)^{-s}f\|^2_{L^2}+(2M)^2\|(|x|+1)^{s}(P(h)-E\pm i\theta)f\|^2_{L^2}.
\end{equation}
We rewrite (\ref{eq:4.2}) in the form
\begin{equation}\label{eq:4.3}
M\theta^{1/2}\|f\|_{L^2}\le \frac{1}{2}\|(|x|+1)^{-s}f\|_{L^2}+
2M^2\|(|x|+1)^{s}(P(h)-E\pm i\theta)f\|_{L^2}.
\end{equation}
We now combine (\ref{eq:4.1}) and (\ref{eq:4.3}) to get
\begin{equation}\label{eq:4.4}
\|(|x|+1)^{-s}f\|_{L^2}\le 4M^2\|(|x|+1)^{s}(P(h)-E\pm i\theta)f\|_{L^2}.
\end{equation}
It follows from (\ref{eq:4.4}) that the resolvent estimate
\begin{equation}\label{eq:4.5}
\left\|(|x|+1)^{-s}(P(h)-E\pm i\theta)^{-1}(|x|+1)^{-s}
\right\|_{L^2\to L^2}\le 4M^2
\end{equation}
holds for all $0<h\ll 1$, $0<\theta\le \epsilon ha^{-2k}$ and $s$ satisfying (\ref{eq:2.1}).
On the other hand, for $\theta\ge \epsilon ha^{-2k}$ the estimate (\ref{eq:4.5}) holds in a trivial way. Indeed, in this case, 
since the operator $P(h)$ is symmetric, the norm of the resolvent is upper bounded by $\theta^{-1}=
{\cal O}(h^{-2km-2})$. Finally, observe that if (\ref{eq:4.5}) holds for $s$ satisfying (\ref{eq:2.1}),
it holds for all $s>1/2$ independent of $h$. Indeed, given an arbitrary $s'>1/2$ independent of $h$, we can arrange by taking $h$ small
enough that $s$ defined by (\ref{eq:2.1}) is less than $s'$. Therefore the bound (\ref{eq:4.5}) holds with $s$ replaced by $s'$ as desired.

\end{document}